\documentclass[12pt,fleqn]{article}

\usepackage{graphicx}
\usepackage{amssymb}

\newtheorem{theorem}{Theorem}[section]
\newtheorem{proposition}[theorem]{Proposition}

\newtheorem{problem}[theorem]{Problem}
\newtheorem{lemma}[theorem]{Lemma}

\newtheorem{corollary}[theorem]{Corollary}

\newtheorem{conjecture}[theorem]{Conjecture}

\begin{document}

\title{Edge-colorings and circular flow numbers on regular graphs}

\vspace{3cm}
      
\author{  Eckhard Steffen\thanks{
		Paderborn Institute for Advanced Studies in 
		Computer Science and Engineering, 
		Paderborn University,
		Germany;			  
		es@upb.de}}      

\date{}

\maketitle

\begin{abstract}
{\small{The paper characterizes $(2t+1)$-regular graphs with circular flow number $2 + \frac{2}{2t-1}$. For $t=1$ this is 
Tutte's characterization of cubic graphs with flow number 4. The class of cubic graphs is the
only class of odd regular graphs where a flow number separates the class 1 graphs from the class 2 graphs. 
We finally state some conjectures and relate them to existing flow-conjectures.
}}
\end{abstract}  

\section{Introduction}

We consider finite (multi-) graphs $G$ with vertex set $V(G)$ and edge set $E(G)$. 
The set of edges which are incident to vertex $v$ is denoted by $E(v)$. 

Vizing \cite{Vizing_1965} proved that the edge-chromatic number $\chi'(G)$ of a graph $G$ with maximum vertex degree $\Delta(G)$ is an 
element of $\{\Delta(G), \dots, \Delta(G) + \mu(G)\}$, where $\mu(G)$ is the maximum multiplicity of an edge of $G$. We say that $G$ is
a class 1 graph if $\chi'(G) = \Delta(G)$ and it is a class 2 graph if $\chi'(G) > \Delta(G)$.

An orientation $D$ of $G$ is an assignment of a direction to each edge, and for $v \in V(G)$, 
$E^-(v)$  is the set of edges of $E(v)$ with head $v$ and $E^+(v)$ is the set of edges with tail $v$.
The oriented graph is denoted by $D(G)$.

A nowhere-zero $r$-flow $(D(G),\phi)$ on $G$ is an orientation $D$ of $G$ together with a function $\phi$
from the edge set of $G$ into the real numbers such that 
$1 \leq |\phi(e)| \leq r-1$, for all $e \in E(G)$, and 
$\sum_{e \in E^+(v)}\phi(e) = \sum_{e \in E^-(v)}\phi(e), \textrm{ for all } v \in V(G)$.
If we reverse the orientation of an edge $e$ with and replace the flow value by $-\phi(e)$, then we
obtain another nowhere-zero $r$-flow on $G$. Hence if there exist an orientation of the edges of $G$ 
such that $G$ has a nowhere-zero $r$-flow, then $G$ has a nowhere-zero $r$-flow for any orientation. 
Thus the question for which values $r$ a graph has a nowhere-zero $r$-flow is a question about graphs,
not directed graphs. Furthermore, $G$ has always an orientation such that all flow values are positive.  
The circular flow number of $G$ is $\inf\{ r | $G$ \mbox{ has a nowhere-zero $r$-flow} \}$, and it is denoted by $F_c(G)$.
It is known, that $F_c(G)$ is always a minimum and that it is a rational number.  

If $G$ has a nowhere-zero flow, then it is 
bridgeless. Tutte \cite{Tutte_1954} conjectured that this necessary structural requirement is also a sufficient condition for a graph to have a 
nowhere-zero 5-flow.  It is easy to see that this conjecture is equivalent to its restriction on cubic graphs. For $i \in \{3,4\}$ there
are characterizations of cubic graphs with nowhere-zero $i$-flow. These results are due to Tutte \cite{Tutte_1949}\cite{Tutte_1954}, see also
\cite{Minty_1967}.

\begin{theorem} [ \cite{Tutte_1949}\cite{Tutte_1954}] \label{Tutte_character}
1) A cubic graph $G$ is bipartite if and only if $F_c(G) = 3$.\\
2) A cubic graph $G$ is a class 1 graph  if and only if $F_c(G) \leq 4$.
\end{theorem}

The following theorem generalizes Theorem \ref{Tutte_character}.1 to $(2t+1)$-regular graphs.

\begin{theorem} [\cite{Steffen_2001}] \label{gap_thm} Let $t \geq 1$ be an integer. 
A $(2t+1)$-regular graph $G$ is bipartite if and only if $F_c(G) = 2 + \frac{1}{t}$. Furthermore,
if $G$ is not bipartite, then $F_c(G) \geq 2 + \frac{2}{2t-1}$.
\end{theorem}

Flow numbers of graphs have attracted considerable attention over the last decades. 
Pan and Zhu \cite{Pan_Zhu_2003} proved that for every rational number $r$ with
$2 \leq r \leq 5$ there is a graph $G$ with $F_c(G) = r$. This result is used in \cite{Schubert_Steffen_2012} to prove the following theorem.

\begin{theorem} [\cite{Schubert_Steffen_2012}] \label{ms_es_2012}
For every integer $t\geq1$ and every rational number $r \in \{2+\frac{1}{t}\} \cup [2+\frac{2}{2t-1};5]$,
there exists a $(2t+1)$-regular graph $G$ with $F_{c}(G)=r$.
\end{theorem}

If $G$ is a cubic graph then $F_c(G) \leq 4$ if and only if $G$ is class 1. Hence, Theorem \ref{Tutte_character}.2 
implies that the flow number 4  separates class 1 and class 2 cubic graphs from each other. 
This paper generalizes Theorem \ref{Tutte_character}.2 to $(2t+1)$-regular graphs. 
We further show that the case of cubic graphs is exceptional in the sense that for every $t > 1$ there is no
flow number that separates $(2t+1)$-regular class 1 graphs and class 2 graphs. However, our results imply that
a $(2t+1)$-regular graph $G$ with $F_c(G) \leq 2+\frac{2}{2t-1}$ is a class 1 graph. We further conjecture that a 
$(2t+1)$-regular graph $H$ with $F_c(H) > 2 + \frac{2}{t}$ is a class 2 graph. We relate this conjecture to other conjectures 
on flows on graphs.

\section{A characterization of $(2t+1)$-regular graphs with circular flow number $\leq 2 + \frac{2}{2t-1}$}

For the proofs of the following results we will use
the concept of  balanced valuations which was introduced by Bondy \cite{Bondy} and Jaeger \cite{Jaeger_75}.
A balanced valuation of a graph
$G$ is a function $w$ from $V(G)$ into the real numbers such that
for all $X \subseteq V(G)$: $| \sum_{v \in X} w(v) | \leq | \partial_G(X) |$, where 
$\partial_G(X)$ is the set of edges with precisely one end in $X$. For $v \in V(G)$
let $d_G(v)$ be the degree of $v$ in the undirected graph $G$.
The following theorem relates integer flows to balanced valuations.

\begin{theorem} [\cite{Jaeger_75}] \label{Thm_Jaeger_75} 
Let $G$ be a graph with orientation $D$ and $r > 2$. Then $G$ 
has a nowhere-zero $r$-flow $(D(G),\varphi)$ if and only if there 
is a balanced valuation $w$ of $G$ such that for all $v \in V(G)$
there is an integer $k_v$ such that $k_v \equiv d_G(v) \bmod 2$ and 
$w(v) = k_v \frac{r}{r-2}$.
\end{theorem}

Furthermore, we need the following result (Theorem 1.1 in \cite{Steffen_2001}).

\begin{lemma} [\cite{Steffen_2001}] \label{n/k}
	Let $n,k$ be integers such that $1 \leq k \leq n$. A graph $G$ has a nowhere-zero $(1 + \frac{n}{k})$-flow
	if and only if $G$ has a nowhere-zero $(1 + \frac{n}{k})$-flow $\phi$ such that for each $e \in E(G)$ there is an integer $m$ such that
		$\phi(e) = \frac{m}{k}$.
\end{lemma}

Note that a cubic graph $G$ is 3-edge-colorable if and only if it has a 1-factor $F$ such that $G-F$ is bipartite. 

\begin{theorem}
Let $t \geq 1$ be an integer. A non-bipartite $(2t+1)$-regular graph $G$ has a 1-factor $F$ such that $G-F$ is bipartite 
if and only if $F_c(G) = 2 + \frac{2}{2t-1}$.
\end{theorem}
{\bf Proof.} ($\leftarrow$) Let $F_c(G) = 2 + \frac{2}{2t-1}$. By Lemma \ref{n/k} there is a $(2 + \frac{2}{2t-1})$-flow $\phi$ with
$\phi(e) \in \{1, 1 + \frac{1}{2t-1}, 1 + \frac{2}{2t-1}\}$ for each $e \in E(G)$. 
Let $F = \{e : \phi(e) = 1 + \frac{1}{2t-1}\}$. We claim that $F$ is a 1-factor of $G$ and
$G-F$ is bipartite. Let $v \in V(G)$ and $|E^+(v)| > |E^-(v)|$. 

Suppose (to the contrary) that 
$\sum_{e \in E^+(v)} \phi(e) > t+1 + \frac{1}{2t-1}$. Then there is an edge $e' \in E^-(v)$ such that
$\phi(e') > \frac{1}{t}(t+1+ \frac{1}{2t-1}) = 1 + \frac{1}{t} + \frac{1}{t(2t-1)} = 1 + \frac{2}{2t-1}$, a contradiction.
Hence, $\sum_{e \in E^+(v)} \phi(e) \leq t+1 + \frac{1}{2t-1}$, $|E^+(v)| = t+1 = |E^-(v)| + 1$.

Furthermore, $|E^+(v) \cap F| \leq 1$. We show that if $|E^+(v) \cap F| = 1$, then $|E^-(v) \cap F| = 0$.
If there is an edge in $E^+(v) \cap F$,
then $t+1 + \frac{1}{2t-1} = \sum_{e \in E^+(v)} \phi(e) = \sum_{e \in E^-(v)} \phi(e)
	\leq t(1 + \frac{2}{2t-1})  = t+1 + \frac{1}{2t-1}$. Hence all edges of $E^-(v)$ have flow value $1 + \frac{2}{2t-1}$,
and $|E^-(v) \cap F| = 0$.

Next we show that if $|E^+(v) \cap F| = 0$, then $|E^-(v) \cap F| = 1$.
If $|E^+(v) \cap F| = 0$, then all edges of $E^+(v)$ have flow value 1. Hence
there are non-negative integers $t_1$, $t_2$, $t_3$ such that $t_1 + t_2 + t_3 = t$ and 
 $t+1 = \sum_{e \in E^+(v)} \phi(e) = \sum_{e \in E^-(v)} \phi(e) = t_1  + t_2 (1 +  \frac{1}{2t-1}) + t_3 (1 +  \frac{2}{2t-1})
	= t + \frac{t_2}{2t-1} + \frac{2t_3}{2t-1}$. Hence, $\frac{t_2}{2t-1} + \frac{2t_3}{2t-1} = 1$ which is 
equivalent to $2t_1 + t_2 = 1$. Thus, $t_1 = 0$, $t_2 = 1$ and therefore,  $|E^-(v) \cap F| = 1$.

It remains to show that $E(v) \cap F \not = \emptyset$. But if $E(v) \cap F = \emptyset$, then $|E^+(v) \cap F| = 0$ and
therefore, $|E^-(v) \cap F| = 1$. Thus $E(v) \cap F \not = \emptyset$, a contradiction. Hence $F$ is a 1-factor of $G$. 

The orientation of the edges induces a 2-coloring of $V(G)$. Let $x$ be a black vertex if $|E^+(x)| = t+1$ and let it be a white vertex
if $|E^+(x)| = t$. 

Let $e \in E(G) - F$ be an edge which is incident to the vertices $v$, $w$, and assume that $e \in E^+(v) \cap E^-(w)$.
We will show that $v$ and $w$ receive different colors. Note that $\phi(e) \in \{1, 1+ \frac{2}{2t-1}\}$.

Suppose to the contrary that $v$ and $w$ have the same color, say both are colored black. Then $|E^+(w)| = t+1$.
If $\phi(e) = 1$ then - since $e \in E^-(w)$ - 
it follows that $\sum_{e \in E^-(w)} \phi(e) \leq 1 + (t-1) (1 +  \frac{2}{2t-1}) < t+1 \leq \sum_{e \in E^+(w)} \phi(e)$, a contradiction.
If $\phi(e) = 1 + \frac{2}{2t-1}$, then - since $e \in E^+(v)$ - it follows that 
$\sum_{e \in E^+(v)} \phi(e) > t+1 + \frac{1}{2t-1}$, a contradiction. 

If both vertices $v$ and $w$ are white, then we deduce a contradiction analogously. Hence, the two vertices of any edge of $G-F$
are in different color classes. Thus, $G-F$ is bipartite. \\[.2cm]
($\rightarrow$) If $G-F$ is a bipartite $2t$-regular graph, then $V(G)$ can be partitioned into two sets $A$ and $B$ with
$|A| =|B|$ and every edge of $G-F$ is incident to one vertex of $A$ and to one vertex of $B$.
Let $w(v) = 2t$ if $v \in A$ and $w(v) = -2k$ if $v \in B$. We claim that $w$ is a balanced valuation on $G$. Let $X \subseteq V(G)$,
$X \cap A = X_A$, $X \cap B = X_B$, and $|X_A|=a$, $|X_B|=b$. We assume that $a \geq b$. It holds that 
$|\partial_G(X)| \geq 2t(a - b) = |\sum_{v \in X} w(v) |$. Hence $G$ has a nowhere-zero $(2+\frac{2}{2t-1})$-flow by Theorem \ref{Thm_Jaeger_75}.
Since $G$ is not bipartite it follows with Theorem \ref{gap_thm} that $F_c(G) = 2+\frac{2}{2t-1}$.
\hfill $\square$

\begin{corollary}
Let $t \geq 1$ be an integer. A $(2t+1)$-regular graph $G$ has a nowhere-zero $(2 + \frac{2}{2t-1})$-flow
if and only if $G$ has a 1-factor $F$ such that $G-F$ is bipartite.
\end{corollary}

\begin{corollary} \label{flow_edge_coloring}
Let $t \geq 1$ be an integer and $G$ be $(2t+1)$-regular graph. If  $F_c(G) \leq 2 + \frac{2}{2t-1}$, then $G$ is a class 1 graph. 
\end{corollary}

\section{Circular flow numbers of class 2 graphs}

Corollary \ref{flow_edge_coloring} generalizes only one direction of Theorem \ref{Tutte_character}.2. The other
direction is already false for $t \geq 2$. In \cite{Steffen_2001} it is shown that $F_c(K_{2t+2}) = 2 + \frac{2}{t}$ for the
complete graph $K_{2t+2}$ on $2t+2$ vertices. Hence, for each $t \geq 2$, there are $(2t+1)$-regular class 1 graphs 
whose circular flow number is greater than $2 + \frac{2}{2t-1}$.

\begin{proposition} \label{prop}
For every integer $t > 1$ and every rational number $r \in \{2+\frac{1}{t-1}\} \cup [2+\frac{2}{2t-3};5]$,
there exists a $(2t+1)$-regular class 2 graph $G$ with $F_{c}(G)=r$.
\end{proposition}

{\bf Proof.} Let $t>1$. By Theorem \ref{ms_es_2012}, for every $r \in \{2+\frac{1}{t-1}\} \cup [2+\frac{2}{2t-3};5]$ there
is a $(2t-1)$-regular graph $G_r$ with $F_c(G_r) = r$. Fix $G_r$ and let $V(G_r) = \{v_1, \dots , v_n\}$.
Let $K_2^{2t+1}$ be the graph on two vertices $u$ and $v$ which are connected by $2t+1$ edges. Let $H_{2t+1}$
be the graph which is obtained from $K_2^{2t+1}$ by subdividing an edge by a vertex $x$. 
For $i \in \{1, \dots,n\}$ let $H_{2t+1}^i$ be a copy of $H_{2t+1}$ with bivalent vertex $x_i$. 
For $t > 1$ let $G'_r$ be the $(2t+1)$-regular graph which is obtained from $G_r$ and $H_{2t+1}^1, \dots, H_{2t+1}^n$
by identifying the vertices $v_i$ of $G_r$ and $x_i$ of $H_{2t+1}^i$ for each $i \in \{1, \dots, n\}$. 
Since $G'_r$ has an odd edge-cut of cardinality smaller than $2t+1$ 
it follows that $G_r'$ is a class 2 graph. Furthermore, $F_c(G'_r) = r$. 
\hfill $\square$

\begin{proposition} \label{prop_2}
For every integer $t > 1$ there are $(2t+1)$-regular graphs $G_1$ and $G_2$ such that $G_1$ is a class 1 graph, 
$G_2$ is a class 2 graph, and $F(G_1) = F(G_2) = 2+\frac{2}{k}$.
\end{proposition}

{\bf Proof.} Let $t > 1$ and $G_1 = K_{2t+2}$. For
$t = 2$ we have $2+\frac{1}{t-1} = 3$ and for $t \geq 3$ holds $2 + \frac{2}{2t-3} \leq 2 + \frac{2}{t}$. Hence, the statement
follows with Proposition \ref{prop}. \hfill $\square$

A $(2t+1)$-regular graph $G$ is a $(2t+1)$-graph if $|\partial_G(X)| \geq 2t+1$ for every $X \subseteq V(G)$ with $|X|$ is odd.
If
$F_c(G) < 2 + \frac{1}{t-1}$, then $G$ must be a $(2 + \frac{1}{t})$-graph. We show that such graphs exist.

Let $G$ be a graph, $F \subseteq E(G)$,
and $F'$ be a copy of $F$. We say that $G'$ is the graph obtained from $G$ by adding $F$ if $V(G')=V(G)$, and $E(G') = E(G) \cup F'$. Let $P$
denote the Petersen graph. The following result is a simple consequence of Theorem 3.1 in \cite{Stefan_es_1999}.

\begin{lemma} [\cite{Stefan_es_1999}] \label{Stefan_es}
Let $k \geq 0$ be an integer. If $G$ is a $(k+3)$-regular graph obtained from $P$ by adding $k$ 1-factors of $P$, then $G$ is class 2.
\end{lemma}

Note, that the graphs of Lemma \ref{Stefan_es} are $(k+3)$-graphs.

\begin{figure}\label{Petersen_fig} 
\centering
	\includegraphics[height=5cm]{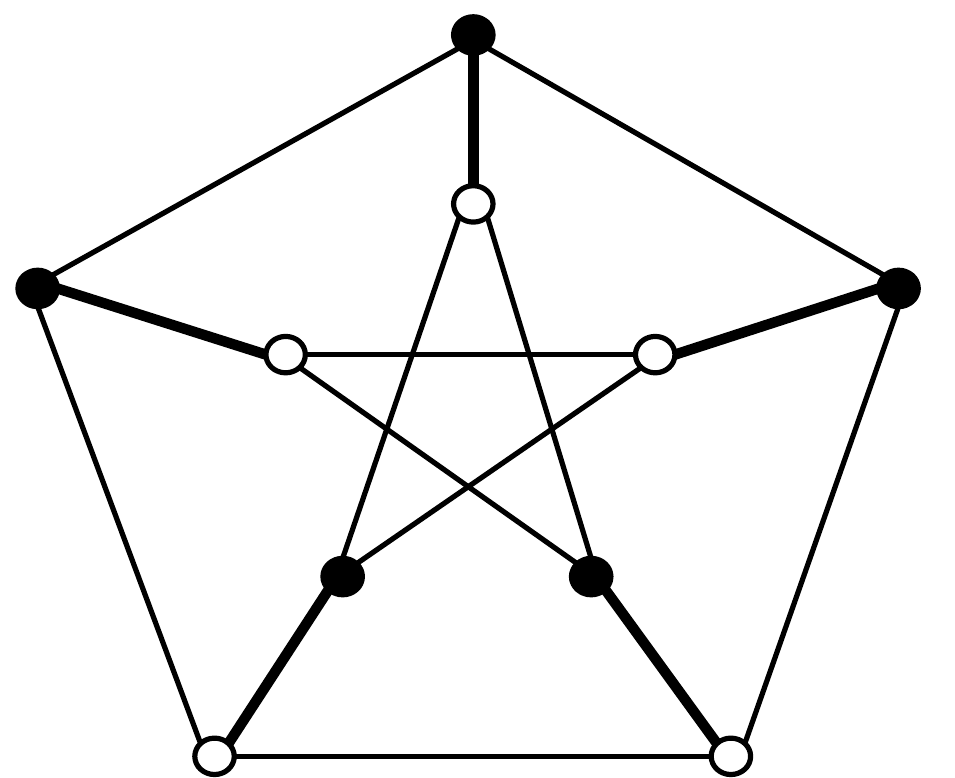}
\caption{The Petersen graph with a vertex 2-coloring. \label{Petersen_fig} }
\end{figure}

\begin{theorem} \label{Petersen_first_step}
For every integer $t \geq 1$ there is a $(2t+1)$-graph $G$ which is a class 2 graph and $F_c(G) = 2 + \frac{3}{3t-2}$.
\end{theorem}

{\bf Proof.} It is well known that $F_c(P)=5$, c.f.~\cite{Steffen_2001}. 
Let the vertices of $P$ be labeled black and white as shown in Figure \ref{Petersen_fig}.
Let $A$ be the set of white vertices and $B$ be the set of black vertices.
It is easy to verify that $w(v) = \frac{5}{3}$ if $v$ is white and $w(v) = - \frac{5}{3}$ if $v$ is black is a balanced valuation on $P$
which corresponds to a nowhere-zero 5-flow on $P$ by Theorem \ref{Thm_Jaeger_75}. Let $F$ be the 1-factor of $P$ which is indicated by the 
bold edges in Figure \ref{Petersen_fig}. Note that if $e \in F$ and $e=xy$, then $x \in A$ if and only if $y \in B$. Let $P_{2t+1}$ be the
$(2t+1)$-graph which is obtained from $P$ by adding $(2t-2)$ copies of $F$. By Lemma \ref{Stefan_es}, $P_{2t+1}$ is a
class 2 graph.

Let $X \subseteq V(P_{2t+1})$, $|\partial_{P_{2t+1}}(X) \cap F| = d$, and $|A \cap X| = a$, $|B \cap X| = b$. 
We assume that $a \geq b$. Since any two vertices of an edge of $F$ belong to different classes it 
follows that $a-b \leq d$. Hence, 
$|\partial_{P_{2t+1}}(X)|  \geq  (2t-2)d + |\partial_P(X)|
	\geq (2t-2)(a-b) + \frac{5}{3}(a-b)
 	\geq (2t - \frac{1}{3})(a-b)$.

Thus, $w_t$ with $w_t(v) = 2t - \frac{1}{3}$ if $v \in A$ and $w_t(v) = -(2t - \frac{1}{3})$ if $v \in B$ is
a balanced valuation on $P_{2t+1}$. Since every partition of $V(P)$ into two classes of cardinality 5 has one class 
which induces a connected component with at least three vertices, it follows that there is no balanced valuation $w'$ on $P_{2t+1}$ with
$|w'(v)| > |w(v)|$. Hence, $F_c(P_{2t+1}) = 2 + \frac{3}{3t-2}$ by Theorem \ref{Thm_Jaeger_75}.   
\hfill $\square$

The results show that for every $t > 1$ there is no flow number that separates $(2t+1)$-regular class 1 graphs from class 2 graphs.
For an integer $t \geq 1$ let 
\[  \Phi(2t+1) = \inf \{F_c(G) : G \textrm{ is a } (2t+1)\textrm{-regular class 2 graph} \}. \]

\begin{corollary} \label{Petersen_first_step}
For every integer $t \geq 1$$:$ $\Phi(2t+1) \leq 2 + \frac{3}{3t-2}$.
\end{corollary}

For cubic graphs $(t=1)$ we have $\Phi(3) = 4$ $(= \frac{2}{2t-1})$. We think that this bound is the right one, and
that the bound of Corollary \ref{flow_edge_coloring} cannot be improved.

\begin{conjecture} \label{inf}
For every integer $t \geq 1$$:$ $\Phi(2t+1) = 2 + \frac{2}{2t-1}$.
\end{conjecture}

The next problem is motivated by Proposition \ref{prop_2}. Furthermore, if it has a positive answer, then
Conjecture \ref{inf} is true.

\begin{problem} \label{problem_1_2}
Is it true that for every integer $t > 1$ and every rational number $r$ with $2 + \frac{2}{2t-1} < r \leq  2 + \frac{2}{t}$
there are $(2t+1)$-regular graphs $H_1$ and $H_2$ such that $H_1$ is class 1, $H_2$ is class 2, and $F_c(H_1) = F_c(H_2) = r$. 
\end{problem}

Let $t \geq 1$ be an integer. 
Corollary \ref{flow_edge_coloring} determines a bound such that all $(2t+1)$-regular graphs with flow number smaller or equal to this bound
are class 1 graphs. We think that there is another flow number such that all $(2t+1)$-regular graphs with flow number greater than this number
are class 2 graphs. 

\begin{conjecture} \label{Conj_upper bound class 1}
Let $t \geq 1$ be an integer and $G$ a $(2t+1)$-regular graph. If $G$ is a class 1 graph, then $F_c(G) \leq 2 + \frac{2}{t}$.
\end{conjecture}

If Conjecture \ref{Conj_upper bound class 1} is true, then the separation of cubic class 1 and class 2 graphs by the flow number 4 is just
due to the fact that $\frac{2}{t} = \frac{2}{2t-1}$ if and only if $t=1$. However, Tutte's 3-flow conjecture is equivalent to the statement that 
$F_c(G) \leq 3$ for every 5-graph $G$. It might be that such a statement is true for each $t >1$.

\begin{conjecture} \label{Conj_upper bound 2t+1-graphs}
Let $t > 1$ be an integer. If $G$ is a $(2t+1)$-graph, then $F_c(G) \leq 2 + \frac{2}{t}$.
\end{conjecture}

Clearly, if Conjecture \ref{Conj_upper bound 2t+1-graphs} is true, then Conjecture \ref{Conj_upper bound class 1} is true. 
Furthermore, if it is true for even $t$, say $t = 2t'$,  then Jaeger's \cite{Jaeger_81} conjecture is true for $(4t'+1)$-regular
graphs. Jaeger \cite{Jaeger_81} conjectured that
every $4t'$-connected graph has a $(2 + \frac{1}{t'})$-flow.

\end{document}